
\input amstex
\documentstyle{amsppt}
\magnification=\magstep1
\def\bc{{\Bbb C}} 
\def\bp{{\Bbb P}}

\newcount\parno \parno=0
\newcount\prono \prono=1
\def\sec{\S\the\parno.-\ \global\prono=1}
\def\etiqueta{\hbox{(\the\parno.\the\prono)}}
\def\finparrafo{\global\advance\parno by1
\vskip 1truecm\ignorespaces}
\def\cita{\ignorespaces\ \the\parno.\the\prono%
\global\advance\prono by 1}

\NoBlackBoxes

\document

\topmatter

\title
On the topology of germs of meromorphic functions
and its applications
\endtitle
\rightheadtext{Topology of meromorphic germs}
\author
S.M. Gusein--Zade, I. Luengo, A. Melle--Hern\'andez
\endauthor

\address
Moscow State University,
Fa\-cul\-ty of Ma\-the\-ma\-tics and Me\-cha\-nics,Mos\-cow, 119899, Russia.
\endaddress
\email
sabir\@mccme.ru
\endemail

\address
Departamento de \'Algebra,
Facultad de Matem\'aticas,
Universidad Complutense.
E-28040 Madrid, Spain.
\endaddress
\email
iluengo\@eucmos.sim.ucm.es\ \ \
amelle\@eucmos.sim.ucm.es
\endemail

\thanks
First author was partially supported by Iberdrola, INTAS--96--0713, RFBR 98--01--00612.
Last two authors were partially supported by DGCYT PB97-0284-C02-01.
\endthanks

\keywords Germs of meromorphic functions, Milnor fibre, atypical values
\endkeywords
\endtopmatter

\document

Germs of meromorphic functions has recently become an object of study in singularity
theory. T. Suwa (\cite{11}) described versal deformations of
meromorphic germs. V.I. Arnold (\cite{1}) classified meromorphic
germs with respect to certain equivalence relations. The authors (\cite{4})
started a study of topological properties of meromorphic germs.
Some applications of the technique developed in \cite{4} were
described in \cite{5} and \cite{6}.

In \cite{4} the authors elaborated notions and technique which could be applied to compute
such invariants of polynomials as Euler characteristics of fibres and
zeta-functions of monodromy transformations associated with a polynomial
(see \cite{5}). Some crucial basic properties of the notions
related to the topology of meromorphic germs were not discussed there. This
has produced some lack of understanding of the general constructions.
The aim of this note is to partially fill in this gap.
At the same time we describe connections with
some previous results and generalizations of them.

A polynomial $P$ in $n+1$ complex variables defines a map $P$
from the affine complex space $\bc^{n+1}$ to the complex line $\bc$.
It is well known that the map $P$
is a $C^\infty$-locally trivial fibration over the complement
to a finite set in the line $\bc$. The smallest of such sets is called
the {\it bifurcation set} or the set of {\it atypical values}
of the polynomial $P$. One is interested in describing the
topology of the fibre of this fibration and its behaviour under
monodromy transformations corresponding to loops around atypical
values of the polynomial $P$. The monodromy transformation
co\-rres\-pon\-ding to a circle of big radius which contains
all atypical values (the monodromy transformation of the polynomial
$P$ at infinity) is of particular interest.

The initial idea was to reduce calculation of the zeta-function of the
monodromy transformation at infinity (and thus of the Euler characteristic of the
generic fibre) of the polynomial $P$ to local problems associated to different
points at infinity, i.e., at the infinite hyperplane $\bc\bp^n_\infty$ in the
projective compactification $\bc\bp^{n+1}$ of the affine space $\bc^{n+1}.$
The
possibility of such a localization for holomorphic germs was used
in \cite{3}. This localization
can be expressed in terms of an integral with respect to the Euler
characteristic, a notion introduced by the school of V.A. Rokhlin (\cite{12}).
However the results are not apply directly  to a polynomial
function since at a point of the infinite hyperplane $\bc\bp^n_\infty$ a
polynomial function
defines not a holomorphic but a meromorphic germ.
Thus the idea of reducing calculations of the zeta-function of the
monodromy transformation of a polynomial map at infinity to calculations of
local zeta-functions corresponding to different points of the hyperplane
$\bc\bp^n_\infty$ is fustrated by the lack of such notions as
the Milnor fibre, the monodromy transformation, ... for meromorphic germs.
Therefore it was necessary to define corresponding invariants and to elaborate a
technique for their calculation.

\finparrafo

\head\sec Basic properties
\endhead

A {\it meromorphic germ} at the origin in the complex space $\bc^{n+1}$
is a ratio ${f=\frac{P}{Q}}$ of two holomorphic germs $P$ and $Q$ on
$(\bc^{n+1},0).$ For our purposes
the following equivalence relation is appropiate. Two
meromorphic germs ${f=\frac{P}{Q}}$ and $f'=\frac{P'}{Q'}$ are {\it equal} if and only if
$P'=P\cdot U$ and $Q'=Q\cdot U$ for a holomorphic germ $U$ not equal to
zero at the origin: $U(0)\not=0$.

A meromorphic germ ${f=\frac{P}{Q}}$ defines a map from the complement to the indeterminacy locus
$\{P=Q=0\}$ to the complex projective line ${\bc\bp^1}$. For any point
$c\in{\bc\bp^1}$, there exists $\varepsilon_0>0$ such that for any
$\varepsilon$, $0<\varepsilon\leq\varepsilon_0,$ the map
$$f:B_\varepsilon\setminus\{P=Q=0\}\to {\bc\bp^1}$$ is a $C^\infty$ locally trivial
fibration over a punctured neighbourhood of the point $c$ here
 $B_\varepsilon$ denotes the closed
ball of radius $\varepsilon$ centred at the origin in $\bc^{n+1}$ (see
\cite{4}).

\definition{Definition} The fibre $${\Cal M}_f^c
=\{z\in B_\varepsilon\,:\, f(z)=\frac{P(z)}{Q(z)}=c'\,\}$$
of this fibration (for $c'$ close to $c$ enough) is called the
{\it $c$-Milnor fibre} of the meromorphic germ $f.$
\enddefinition

The Milnor fibre ${\Cal M}_f^c$ is a (non-compact) $n$-dimensional complex manifold
with boundary.

\definition{Definition} The monodromy transformation of this
fibration corresponding to a simple
(small) loop around the value $c$ is called the {\it $c$-monodromy
transformation} of the meromorphic germ $f$.
\enddefinition

\definition{Definition} A value $c\in{\bc\bp^1}$ is called {\it typical} if the map
$f:B_\varepsilon\setminus\{P=Q=0\}\to{\bc\bp^1}$ is a $C^\infty$ locally trivial
fibration over a neighbourhood of the point $c$ (including $c$ itself).
\enddefinition

Notice that, if $c$ is a typical value, the corresponding monodromy
transformation is isotopic to the identity.

\proclaim{Theorem 1} There exists a finite set $\Sigma\subset{\bc\bp^1}$ such that for
all $c\in {\bc\bp^1}\setminus \Sigma$ the $c$-Milnor fibres of $f$ are diffeomorphic
to each other and the $c$-monodromy transformations are trivial (i.e., isotopic
to identity). In particular, the set of atypical values is finite.
\endproclaim

\demo{Proof} A {\it resolution of the germ $f$} is a
modification of the space $(\bc^{n+1},0)$ (i.e., a proper analytic map $\pi:{\Cal X}\to
{\Cal U}$ of a smooth analytic manifold $\Cal X$ onto a neighbourhood $\Cal U$
of the origin in $\bc^{n+1}$, which is an isomorphism outside of a proper
analytic
subspace in $\Cal U$) such that the total transform $\pi^{-1}(H)$ of the
hypersurface $H=\{P=0\}\cup\{Q=0\}$ is a normal crossing divisor at each
point of the manifold ${\Cal X}$. We assume that the map $\pi$
is an isomorphism outside of the hypersurface $H$.

The fact that the preimage
$\pi^{-1}(H)$ is a divisor with normal crossings implies that in a
neighbourhood
of any point of it, there exists a local system of coordinates
$y_0,y_1,\ldots,y_n$
such that the liftings $\widetilde{P}=P\circ \pi$ and $\widetilde{Q}=Q\circ \pi$
of the
functions $P$ and $Q$ to the space $\Cal X$ of the modification are equal to
$u\cdot y_0^{k_0}y_1^{k_1}\cdots y_n^{k_n}$ and $v\cdot
y_0^{l_0}y_1^{l_1}\cdots y_n^{l_n}$
respectively, where $u(0)\not=0$ and $v(0)\not=0$,
$k_i$ and $l_i$ are nonnegative.

\remark{Remark 1} The values $0$ and $\infty$ in the projective line ${\bc\bp^1}$
are used as distinguished points for convenience: to have the usual notion
of a resolution of a function for the numerator and for the denominator.
\endremark

One can make additional blow-ups along intersections of pairs of irreducible
components of the divisor $\pi^{-1}(H)$ so that the lifting
${\widetilde f}=f\circ \pi=\frac{\widetilde{P}}{\widetilde{Q}}$ of the function $f$ can
be defined as a holomorphic map from the manifold ${\Cal X}$ to the
complex projective line ${\bc\bp^1}$. This condition means that
$\widetilde{P}=V\cdot P',\,\widetilde{Q}=V\cdot Q'$ where $V$ is a section
of a line bundle, say ${\Cal L},$ over ${\Cal X},$ $P'$ and $Q'$ are
sections of the line bundle ${\Cal L}^{-1}$ and $P'$ and $Q'$ have no
common zeroes on $\Cal X$. Let ${\widetilde f}'=\frac{P'}{Q'}$.

On each component of the divisor $\pi^{-1}(H)$
and on all the intersections of several of them ${\widetilde
f}'$ defines a map to the projective line ${\bc\bp^1}$. These maps have a finite
number of critical values, say $a_1$, $a_2$, ..., $a_s$.

\remark{Remark 2} If the function ${\widetilde f}'$ is constant on a
component of a finite intersection of the irreducible divisors of 
$\pi^{-1}(H)$, then this
constant value is critical. The value the function ${\widetilde f}'$
on an intersection of $n+1$ components (this intersection is zero-dimensional)
should also be considered as a critical value. \endremark

Let $c\in {\bc\bp^1}$ be different from $a_1$, $a_2$, ..., $a_s$. We shall show that
for all $c'$ from a neighbourhood of the point $c$ (including $c$ itself)
the $c'$-Milnor fibres of the meromorphic function $f$ are diffeomorphic
to each other and the $c'$-monodromy transformations are trivial.

Let $r^2(z)$ be the square of the distance from the origin in the space
$\bc^{n+1}$ and let ${\widetilde r}^2(x)=r^2(\pi(x))$ be the lifting of this
function to the space $\Cal X$ of the modification. In order to define the
$c'$-Milnor fibre one has to choose $\varepsilon_0>0 $ (the Milnor radius)
small enough so that the level manifold $\{{\widetilde r}^2(x)=
\varepsilon^2\}$ is transversal to $\{{\widetilde f}'(x)=c'\}$ for all
$\varepsilon$ such that $0<\varepsilon\leq \varepsilon_0$. Let
$\varepsilon_0=\varepsilon_0(c) $ be the Milnor radius for the value $c$.
Since $\{{\widetilde f}'(x)=c\}$ is transversal to
components of the divisor $\pi^{-1}(H)$ and to all their intersections,
then also $\varepsilon_0$ is the Milnor radius for 
all $c'$ from a neighbourhood
of the point $c\in {\bc\bp^1}$ (and the level manifold
$\{{\widetilde f}'(x)=c'\}$ is transversal to components of the divisor
$\pi^{-1}(H)$ and to its intersections). This implies that for such $c'$
the $c'$-Milnor fibres of the meromorphic germ $f$ are diffeomorphic to each
other and the $c'$-monodromy transformations are trivial.
$\square$\enddemo

\remark{Remark 3} The $c$-Milnor fibre for a generic value
$c\in {\bc\bp^1}$ can be called {\it the generic Milnor fibre} of the meromorphic
germ $f$. One can easily see that the generic Milnor fibre of a meromorphic
germ can be considered as embedded into the $c$-Milnor fibre for any value
$c\in{\bc\bp^1}$.
Moreover the Euler characteristic of the generic Milnor fibre of a meromorphic
germ is equal to zero and the zeta-function of the corresponding
monodromy transformation (see \cite{4}) is equal to $(1-t)^0=1.$
\endremark

\finparrafo

\head\sec Isolated singularities
and Euler characteristic of the $0$-Milnor fibre
\endhead

Let $P$ be a polynomial in $n+1$ complex variables. Suppose that the
closure ${\overline V}_{t_0}\subset\bc\bp^{n+1}$ of the level set
$V_{t_0}=\{P=t_0\}\subset\bc^{n+1}$ in the complex projective
space $\bc\bp^{n+1}\supset \bc^{n+1}$ has only isolated singular points. Let
$A_1$, \dots, $A_r$ be those of them which lie in the affine space $\bc^{n+1}$,
and let $B_1$, \dots, $B_s$ be those which lie in the infinite hyperplane
$\bc\bp^n_\infty$. For
$t$ close enough to $t_0$ (and thus generic),
the closure ${\overline V}_{t}\subset\bc\bp^{n+1}$ of the level set
$V_{t}=\{P=t\}\subset\bc^{n+1}$ has no singular points in 
the space $\bc^{n+1}$
and may have isolated singularities only at the points $B_1$, \dots, $B_s$. It is
known that $$\chi(V_t)-\chi(V_{t_0})=(-1)^{n+1}\left(\sum_{i=1}^r
\mu_{A_i}(V_{t_0}) +\sum_{j=1}^s(\mu_{B_j}({\overline
V}_{t_0})-\mu_{B_j}({\overline V}_{t}))\right).\qquad (1)$$

We shall formulate a somewhat more general statement about meromorphic germs.
This statement together with the formula for the difference
of the Euler characteristics of the
generic level set of a polynomial $P$ and of a special one in terms of
meromorphic germs defined by the polynomial $P$ (see
\cite{6}) gives $(1).$

\proclaim{Theorem 2} Let ${f=\frac{P}{Q}}$ be a germ of meromorphic function on the space
$(\bc^{n+1},0)$ such that the numerator $P$ has an isolated critical point
at the origin and, if $n=1$, the germs of the curves $\{P=0\}$ and $\{Q=0\}$
have no common irreducible components. Then, for a generic $t\in\bc,$
$$\chi({\Cal M}_{f}^0)=(-1)^n(\mu(P,0)-\mu(P+tQ,0)).$$
 Here
$\mu(g,0)$ stands for the usual Milnor number of the holomorphic germ $g$
at the origin. \endproclaim

\demo{Proof} The Milnor fibre ${\Cal M}_{f}^0$ of the meromorphic germ
$f$ has the following description. Let $\varepsilon$ be small enough (and
thus a Milnor radius for the holomorphic germ $P$). Then
$${\Cal M}_{f}^0=B_\varepsilon(0)\cap(\{P+tQ=0\}\setminus\{P=Q=0\}),$$ for
$t\not=0$ with $|t|$ small enough (and thus $t$ generic). Note that the
zero-level set $\{P+tQ=0\}$ is non-singular outside of the origin for
$|t|$ small enough. The space $B_\varepsilon(0)\cap\{P=Q=0\}$ is homeomorphic
to a cone and therefore its Euler characteristic is equal to $1$. Therefore
$$\chi({\Cal M}_{f}^0)=\chi(B_\varepsilon(0)\cap(\{P+tQ=0\})-1.$$ Now
Theorem 2 is a consequence of the following well known fact
(see, e.g., \cite{2}).

\noindent {\bf Statement.} Let $P:(\bc^{n+1},0)\to(\bc,0)$ be a germ of a
holomorphic function with an isolated critical point at the
origin and let $P_t$ be any deformation of $P$ ($P_0=P$). Let
$\varepsilon$ be small enough. Then for $|t|$ small enough
$$(-1)^n(\chi(B_\varepsilon(0)\cap\{P_t=0\})-1)$$ is equal to the number of
critical points of $P$ (counted with multiplicities) which split from the zero
level set, i.e., to $$\mu(P,0)-\sum_{Q\in\{P_t=0\}\cap B_\varepsilon}
\mu(P_t,Q).\qquad\square$$
\enddemo

\example{Example 1} The following example shows the necessity of the condition
that, for $n=2,$ the curves $\{P=0\}$ and $\{Q=0\}$ have no common components.
Take $P=xy$ and $Q=x$. \endexample

\example{Example 2} In Theorem 2 the difference of Milnor numbers
(up to a sign) may be replace by (the equal) difference of Euler
characteristics of the corresponding Milnor fibres (of the germs $P$ and
$P+tQ$). However the formula obtained this way is not correct if the
germ $P$ has a nonisolated critical point at the origin. It is shown
by the example $f=\frac{x^2+z^2y}{z^2}$. \endexample

The formula $(1)$ is a direct consequence of Theorem 2 and the formula $(2)$
of Theorem 2 from \cite{6}.

\finparrafo

\head\sec Topological triviality of the family $\{P+tQ\}$ and typical values of
meromorphic germs \endhead

Again let ${f=\frac{P}{Q}}$ be a meromorphic germ on $(\bc^{n+1},0)$ such that the
holomorphic germ $P$ has an isolated critical point at the origin.

\proclaim{Theorem 3} The value $0$ is typical for the meromorphic germ $f$ if
and only if $\chi({\Cal M}_{f}^0)=0$. \endproclaim

\demo{Proof} `` Only if\," follows from the definition and Theorem 2.

`` If\," is a consequence of the result of A. Parusi\'nski \cite{8}
(or rather of its proof). He has proved that, if $\mu(P)=\mu(P+tQ)$
for $|t|$ small enough, then the family of maps $P_t=P+tQ$ is
topologically trivial. In particular the family of germs of hypersurfaces
$\{P_t=0\}$ is topologically trivial. For $n\not=2$ this was proved by
L\^e D.T. and C.P. Ramanujam \cite{7}. However in order to apply the
result to the present situation it is necessary to have a topological
trivialization of the family $\{P_t=0\}$ which preserves the subset $\{P=Q=0\}$
and is smooth outside the origin. For the family $P_t=P+tQ$, such a
trivialization was explicitly constructed in \cite{8} without
any restriction on the dimension. $\qquad\qquad\square$\enddemo

\example{Example 3} If the germ of the function $P$ has a non-isolated
critical point at the origin then this characterization is no longer true.
Take, for example, $P(x, y)=x^2y^2$ and $Q(x, y)=x^4+y^4$.
\endexample

\bigbreak
\finparrafo

\head\sec A generalization of the Parusi\'nski--Pragacz formula for the Euler
characteristic of a singular hypersurface
\endhead

Let $X$ be a compact complex manifold and let $\Cal L$ be a holomorphic
line bundle on $X.$ Let $s$ be a section of the bundle $\Cal L$ not
identically equal to zero, and $Z:=\{s=0\}$ is its zero locus (a hypersurface
in the manifold $X$). Let $s'$ be another section of the bundle $\Cal L$
whose zero locus $Z'$ is nonsingular and transversal to a Whitney
stratification of the hypersurface $Z$. A. Parusi\'nski and P. Pragacz
have proved (see \cite{9}, Proposition 7) a statement
which in terms of \cite{6} can be written as follows
$$\chi(Z')-\chi(Z)=\int_{Z\setminus Z'} (\chi_x(Z)-1)\,d\chi,\qquad (2)$$
where $\chi_x(Z)$ is the Euler characteristic of the Milnor fibre of the germ
of the section $s$ at the point $x$ (the definition of the integral with respect
to the Euler characteristic can be found in \cite{12} or \cite{6}).

We shall indicate a more general formula which
includes this as a particular case.

\proclaim{Theorem 4} Let $s$ be as above and let $s'$ be a section of the
bundle $\Cal L$
whose zero locus $Z'$ is non-singular. Let $f$ be the meromorphic
function $s/s'$ on the manifold $X$. Then
$$\chi(Z')-\chi(Z)=\int_{Z\setminus Z'} (\chi_x(Z)-1)\,d\chi\,+
\int_{Z\cap Z'}\chi_{f,x}^0 d\chi,\qquad (3)$$
where $\chi_{f,x}^0$ is the Euler characteristic of the $0$--Milnor
fibre of the meromorphic germ $f$ at the point $x$.

\endproclaim

\demo{Proof} Let $F_t$ be the level set $\{f=t\}$ of the (global) meromorphic
function $f$ on the manifold $X$ (with indeterminacy set $\{s=s'=0\}$), i.e.,
$F_t=\{s-ts'=0\}\setminus \{s=s'=0\}.$ By \cite{6}, for a generic value
$t$ one has
$$\chi(F_{gen})-\chi(F_0)=\int_{F_0} (\chi_{f,x}^0(Z)-1)\,d\chi\,+
\int_{\{s=s'=0\}}\chi_{f,x}^0 d\chi,$$
where $\chi_{f,x}^0$ is the Euler characteristic of the $0$-Milnor fibre of
the meromorphic germ $f$ at the point $x.$ One has $F_0=Z\setminus (Z\cap
Z'),\, F_\infty=Z'\setminus (Z\cap Z')$ and in this case $F_\infty$ is a
generic level set of the meromorphic function $f$ (since its closure is non-singular).
Therefore $\chi(F_0)=\chi(Z)-\chi(Z\cap Z'),$
$\chi(F_{gen})=\chi(Z')-\chi(Z\cap Z').$ Finally, for $x\in F_0$, the
germ of the function $f$
at the point $x$ is holomorphic and thus $\chi_{f,x}^0=\chi_x(Z).\qquad
\square$ \enddemo

If the hypersurface $Z'$ is transversal to all strata of a Whitney
stratification of the hypersurface $Z,$ then, for $x\in Z\cap Z',$ the Euler
characteristic $\chi_{f,x}^0=0$ (Proposition 5.1 from \cite{10})
and therefore the formula $(3)$ reduces to $(2).$

\bigbreak
\finparrafo


\enddocument